\documentclass[11.1pt,fleqn]{article}

\usepackage[english]{babel}
\usepackage{amsmath, amssymb, multicol, amsfonts, amsthm, amscd, epsfig, enumerate}
\usepackage[latin1]{inputenc}
\usepackage[T1]{fontenc}
\usepackage{graphicx}
\usepackage{geometry}
\usepackage{algorithm}
\usepackage{algpseudocode}
\algrenewcommand\algorithmicrequire{\textbf{Initialize:}}

\newtheorem{theorem}{Theorem}

\newtheorem{corollary}[theorem]{Corollary}

\newtheorem{proposition}[theorem]{Proposition}
\newtheorem{remark}[theorem]{Remark}

\begin{document}

\author{
	\emph{\textbf{Alessandro Ramponi}}\\
Economics Department, University of Rome Tor Vergata\\
Via Columbia 2, \ 00133 Rome, \ Italy\\
alessandro.ramponi@uniroma2.it
	\\
	\emph{\textbf{M. Elisabetta Tessitore}}\\
	Economics Department, University of Rome Tor Vergata\\
	Via Columbia 2, \ 00133 Rome, \ Italy\\
	tessitore@economia.uniroma2.it }

\title{Optimal Social and Vaccination Control in the SVIR Epidemic Model}
\date{\today}
\maketitle

\begin{abstract}
In this paper we introduce an approach to the management of infectious disease diffusion through the formulation of a controlled compartmental SVIR (Susceptible-Vaccinated-Infected-Recovered) model. We consider a cost functional encompassing three distinct yet interconnected dimensions: the social cost, the disease cost, and the vaccination cost. The proposed model addresses the pressing need for optimized strategies in disease containment, incorporating both social control measures and vaccination campaigns. Through the utilization of advanced control theory, we identify optimal control strategies that mitigate disease proliferation while considering the inherent trade-offs among social interventions and vaccination efforts. Finally, numerical implementation of the optimally controlled system through the Forward-Backward Sweep algorithm is presented. 
\end{abstract}

\thispagestyle{empty}

\textbf{Keywords: optimal control;  economics; SVIR Epidemic Model}

\textbf{JEL: }

\section{Introduction}

In recent years, numerous disease outbreaks have raised concerns about global public health. Avian influenza, Marburg virus, measles, lassa fever, and other infectious diseases have highlighted the potential dangers they pose. Organizations like the World Health Organization (WHO) have been at the forefront of monitoring and responding to these outbreaks, recognizing the urgent need for effective disease control strategies.

In the wake of the COVID-19 pandemic, the significance of robust response measures has been underscored. The \textit{European Centre for Disease Prevention and Control} has outlined key lessons learned in their technical report titled "Lessons from the COVID-19 pandemic"\footnote{\texttt{https:\slash \slash www.ecdc.europa.eu \slash en\slash publications-data\slash lessons-covid-19-pandemic-may-2023}}. The report identifies four critical components of the response to health threats, emphasizing the pivotal role of analyzing epidemiological data and employing epidemiological modeling to inform decision-making.

Mathematical modeling remains a crucial factor in epidemiology, offering a more profound understanding of the fundamental mechanisms behind the propagation of newly emerging and resurging infectious diseases, as well as proposing efficient strategies for their control. In such a framework, compartmental models, which originated with the Kermack and McKendrick (1927) \cite{KM} Susceptible-Infectious-Recovered (SIR) model and have since undergone extensive development (see e.g. the book Brauer and Castillo-Chavez 2010 \cite{BC}), represent a recognized and established class of dynamic models used to depict the progression of infectious diseases. Within this context, and focusing on the problem of disease containment, we regard on two main tools for its control: the implementation of a set of social restrictions (e.g. lockdown periods) and the deployment of vaccination campaigns. The explicit introduction of a vaccination compartment has been proposed in Liu et al (2008) \cite{XTI}, which consider a modified version of the SIR model, named SVIR,
 in which the new  compartment $V$ is added, where the vaccenees will belong before reaching the immunity and therefore the recovered individuals. 
 
Notably, while the adoption of stringent lockdown measures holds and the availability of a vaccine promise in curbing the spread of the virus, they raise significant concerns regarding their pronounced economic impact, thereby elevating government decision-making to a multifaceted challenge. Optimal control theory of compartmental models offers a solid theoretical framework to capture essential aspects of an optimal disease control policies. By representing the population as distinct compartments and modeling the dynamics of disease transmission, optimal control theory provides a systematic approach to determining the most effective strategies for disease control. In fact, this framework enables the exploration of various control measures and their impact on disease spread, allowing policymakers to optimize interventions based on desired objectives such as minimizing infection rates, reducing economic costs, or maximizing resource utilization.

\medskip

The literature on infectious disease analyzed via optimal control (see e.g. Lenhart and Workman 2007 \cite{LW}) is experiencing rapid and extraordinary development. In such a context, Behncke (2000) \cite{Ben00} represents one of the early endeavors to systematically incorporate a control methodology. In preceding decades, research primarily concentrated on strategies centered around selective isolation and immunization. Abakuks (1973) \cite{A} explored the optimal separation of an infected population under the assumption of instantaneous isolation, while Hethcote and Waltman (1973) \cite{HW73} introduced optimal vaccination strategies in their work. In more recent studies, Ledzewicz and Schattler (2011) \cite{LS} addressed an optimal control problem within the context of a model incorporating both vaccines and treatments in a dynamically expanding population. In a related vein, Gaff and Schaefer (2009) \cite{GS} conducted research encompassing SIR/SEIR/SIRS models, focusing on control parameters that govern vaccination rates, treatments provided to infected individuals, and the potential for quarantine measures to be applied as well. Bolzoni et al. (2017)  \cite{B} conducted an examination of time-optimal control problems concerning the utilization of vaccination, isolation, and culling strategies within the context of a linear framework. In the work of Miclo et al. (2020) \cite{MSW20}, the researchers investigated a deterministic SIR model wherein a social planner exercised control over the transmission rate. This control was aimed at mitigating the transmission rate's natural levels to prevent an undue strain on the healthcare system. Kruse and Strack (2020), as detailed in \cite{KS}, expanded the SIR model to include a parameter subject to the planner's control, influencing disease transmission rates. This parameterization captured political measures such as social distancing and institution and business lockdowns, which, while effective in curtailing disease transmission, often incurred substantial economic and societal costs. This trade-off was modeled by introducing convex cost functions related to the number of infected individuals and reductions in transmission rates. Addressing the complex issue of epidemic management, Federico and Ferrari (2021) \cite{FF} focused on a policymaker's endeavor to curtail the epidemic's spread while concurrently minimizing the associated societal costs within a stochastic extension of the SIR model. Concurrently, Federico et al. (2022) \cite{FFT22} investigated an optimal vaccination strategy utilizing a dynamic programming approach within an SIRS compartmental model. Calvia et al. (2023) \cite{CGLZ22} delved into the control of epidemic diffusion through lockdown policies within an SIRD model, employing a dynamic programming framework for in-depth analysis. In Chen et al. (2022) \cite{CPW22}, a similar compartmental model is employed to investigate the impact of social distancing measures on mitigating COVID-19, examining the situation from both economic and healthcare standpoints. The study utilizes daily pandemic data, including figures for infected, recovered, and deceased individuals, in addition to economic indicators such as mobility and financial market instability indices. The overarching multi-objective is to minimize the risks associated with disease transmission and economic volatility.

The application of an optimal control framework to an SVIR dynamical model has received relatively limited attention within the existing literature. In the work of Ishikawa (2012) \cite{Ish2012}, focus is directed toward a stochastic version of this model, where a thorough analysis of the corresponding stochastic optimal control problem revolves around the vaccination strategy, featuring a quadratic cost function. Witbooi et al. (2015) \cite{WMV} extended the investigation, addressing both deterministic and stochastic optimal control problems within the SVIR model framework. Their approach considers the vaccination rate as a controllable parameter and integrates an additive cost function. Kumar and Srivastava (2017) \cite{KPS} proposed and examined a control problem within this framework, incorporating both vaccination and treatment as control policies. Notably, they introduced a cost function linear in state variables, quadratic in treatment measures, and quartic in vaccination policies, respectively. Similarly, in the study conducted by Garriga et al. (2022) \cite{GMS2022}, the deterministic optimal control problem is explored in the context of a pandemic characterized by two distinct phases. During the initial phase, social restrictions serve as the sole viable containment measures for the disease, while at a subsequent random time, the availability of a vaccine is introduced. Optimal control strategies are thoroughly examined for both phases, encompassing the utilization of one and two control variables, respectively. These analyses also specify the cost function's structural attributes through the incorporation of an utility function.

In this paper we extend the results obtained in \cite{RT23} by adding another dimension to the decision-making process. In our analysis, we postulate a scenario where a disease has already disseminated at an early stage, subsequently followed by the availability of a vaccine. Hence, we employ a SVIR dynamic model to depict the transmission of an infectious disease, which can be influenced by two types of mitigation measures within the purview of a social planner. These measures aim to curtail the rate of contagion within the population to decrease the impact of the disease. The central challenge lies in determining the optimal response by striking a balance between the restrictions that minimize disease prevalence, the vaccination rate and the economic costs associated with the strategies implementation.  To comprehensively account for the impact of these measures, differently from \cite{WMV}, we introduce an explicit cost function that distinctly factors in the expenses associated with handling the infected population, conducting vaccination campaigns, and the economic impact resulting from social restrictions. By specifying the functional form of the cost functional, we are able to characterize the optimal control strategy function through the application of Pontryagin's Maximum Principle, see e.g. Lenhart and Workman (2007) \cite{LW}.

The analysis of the optimally controlled SVIR dynamic, and the corresponding optimal policies, is carried out by implementing the Forward-Backward Sweep algorithm, as described in Lenhart and Workman (2007) \cite{LW}. It is a two-step procedure that proved to be a practical approach to solving a wide range of optimal control problems by iteratively refining the control functions based on the state-costate variables, as established by the necessary conditions given by the maximum principle. The numerical simulations are conducted within two main epidemiological scenarios, characterized by different basic reproduction numbers, corresponding to disease-free and endemic equilibria of the dynamical model, respectively. The main results of our empirical analysis show the consistent reduction in total cost achieved by implementing the optimal policy, compared with the three benchmark strategies considered. 
 
\medskip

The paper is structured as follows. In Section 2, we introduce the SVIR dynamical model and the corresponding control problem, for which we characterize the functional form of the optimal control strategies. The quantitative analysis is reported in Section 3, which summarizes the results of the numerical simulations. In the final Section 4, we highlight our main findings and discuss directions for future research.

\section{Problem formulation} 

\subsection{The SVIR model}\label{basic_svir}

\begin{figure}[h]
\begin{center}
 	\includegraphics[width=5cm, height=3cm]{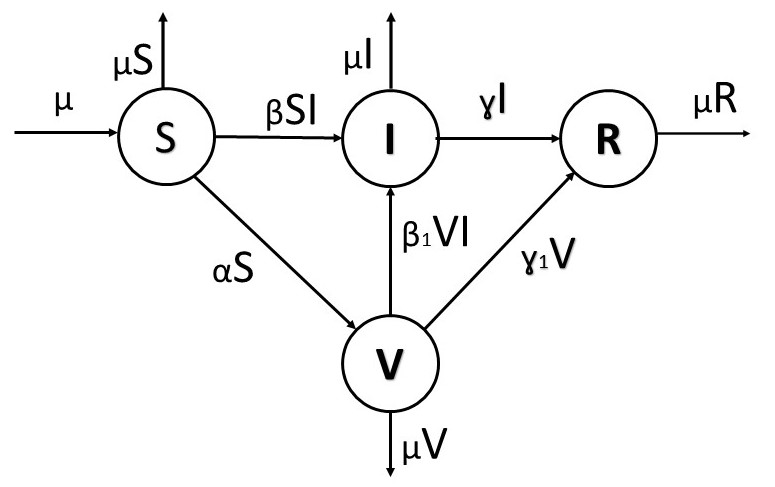}
\end{center}
\caption{The basic SVIR model graph.}
\label{svir_grafo1}
\end{figure}

Liu et al. proposed the SVIR model in \cite{XTI} to extend the well-known SIR model by adding a vaccination program (continuous or impulsive) for the population under study. The model consists of four groups: the Susceptibles S, the Infected I, the Recovered R, and the Vaccinees V, who are those who have started the vaccination process. The fractions of the total population in each group are denoted by S, V, R, and I, respectively.
The model assumes that the disease is transmitted at a rate $\beta $ when the susceptible individuals come into contact with the infected ones, and that the infected individuals recover at a rate $\gamma $. The vaccinated individuals acquire immunity against the disease at a rate $\gamma _1 $  and they can also be infected at a reduced rate $\beta _1 $, which is lower than$\beta $ because some immunity is gained after vaccination. The parameter $\alpha $ represents the rate at which the susceptible individuals join the vaccination program, and $\mu $ is the birth-death rate. Figure \ref{svir_grafo1} illustrates how the population moves among  the four compartments S,V,I,R.

The following system of first order differential equations captures the framework of the continuous vaccination process:
\begin{equation}
\label{svirconst}
\left\{
\begin{array}{ll}
\displaystyle{\frac{dS}{dt}(t)} =-\beta S(t)I(t) -\alpha S(t)+ \mu - \mu S(t) &S(0)=S_0
\\\\
\displaystyle{\frac{dV}{dt}}(t) =\alpha S(t) - \beta_1  V(t)I(t)-\gamma _1V(t)- \mu V(t) &V(0)=V_0
\\\\
\displaystyle{\frac{dI}{dt}}(t) =\beta S(t)I(t) + \beta_1 V(t)I(t)-\gamma I(t) - \mu I(t) &I(0)=I_0
\\\\
\displaystyle{\frac{dR}{dt}}(t)=\gamma _1V(t)+\gamma I(t)- \mu R(t) &R(0)=R_0
\end{array}
\right.
\end{equation}

Let the parameters $\beta, \beta_1, \gamma, \gamma_1, \mu $ be positive real numbers and $\alpha $ be a non-negative real number. We also assume that the initial values $S_0, V_0, I_0, R_0$ are positive real numbers and their sum is equal to 1. These assumptions are made because the model (\ref{svirconst}) describes human populations, and it can be proven that the solutions of the system remain non-negative if the initial values are non-negative, as shown in \cite{XTI}. Furthermore, we note that if we define N(t)=S(t)+V(t)+I(t)+R(t), then we can see from (\ref{svirconst}) that $\frac {dN}{dt}(t) = 0$: therefore $N(t)=N_0 \equiv 1$, for all $t\geq 0$.

The state variable $R$ does not appear in the first three equations of system (\ref{svirconst}), so we can analyze the properties of the system using only the variables $S$, $V$, and $I$. As shown in \cite{XTI}, the SVIR model has a disease free equilibrium $E_0$ (meaning an equilibrium $E_0=(S^*, V^*, I^*)$ such that $I^*\equiv 0$), and an endemic equilibrium $E_+=(S^+, V^+, I^+)$ with $I^+>0$. Furthermore, the \textit{basic reproduction number} determines its long-term behavior.

\begin{equation}\label{R0}
R_0^C = \frac{\mu \beta}{(\mu+\alpha)(\mu+\gamma)}+\frac{\alpha \mu \beta_1 }{(\mu+\gamma_1)(\mu+\alpha)(\mu + \gamma)},
\end{equation}
and it is summarized in the two following Theorems, which are proved in \cite{XTI}:
\begin{theorem}
If $R^C_0 <1$, then the disease free equilibrium $E_0$, which always exists, is locally asymptotically stable  and is unstable if $R^C_0>1$. 
Moeover, $R^C_0 >1$ if and only if system (\ref{svirconst}) has a unique positive equilibrium $E_+$   and it is locally asymptotically stable when it exists.
\end{theorem}

\begin{theorem}
If $R^C_0 \leq 1$, then the disease free equilibrium $E_0$ is globally asymptotically stable. And if $R^C_0 > 1$, the endemic equilibrium $E_+$ is globally asymptotically stable in all the region of feasible model solutions except for the constant solution identically equal to $E_0$.

\end{theorem}


\subsection{The optimal control problem} \label{sect_control}
This Section presents the controlled SVIR model and its related deterministic optimal control problem.

As in standard SVIR models, let S,V,I and R represent  susceptible,  vaccinated ,  infected  and recovered  respectively. By S,V,I and R we denote the percentage of the total population belonging to each group.

We consider a vector control variable $u(\cdot )=(u_0(\cdot ),u_1(\cdot ))$. The component $u_0(\cdot )$  is meant to govern the social restrictions imposed by the social planner on a population until a specific time $T$, which is the final time of government restriction, while the component $u_1(\cdot )$  is meant to govern the rate at which susceptible people are moving into the vaccinees compartment via the vaccination program. The birth and the death rate is assume to be the same and it is denoted by $\mu $.

The control variable $u=(u_{0},u_{1})$ belongs to the admissible set $\cal U$ defined as 
$$
{\cal {U}}=\left \{u:[0,T]\to [0,\overline{u_0}]\times [0,\overline{u_1}] :\mbox{Lebesgue measurable}, \;  \overline{u_0}, \overline{u_1}\in (0,1] \right \}
$$

The control variable $u_0$ allows to \textit{adjust} the rate of transmission of the disease. It captures  the restrictions social planner imposes to govern the speed at which the infection spreads, while  $\varepsilon$ quantifies the vaccine  
ineffectiveness (if $\varepsilon \equiv 0$ no vaccinated gets infected). The goal is to create a scenario where the infection rate $\beta$ is high without control ($u_0=0$), and low with increasing controls. The function $\beta (u_0)$ reflects both the infectiousness of the disease and the social planner's measures to control the infection spread.
 
The control variable $u_1$ captures the cost of the vaccination, which we assume to be proportional to $u_1S$, being the flux of individuals from $S$ to $V$ or, equivalently, the number of new vaccinees individuals in the unit of time. Furthermore we set $\beta_1= \varepsilon \beta$, where $\varepsilon$ quantifies the vaccine ineffectiveness (if $\varepsilon \equiv 0$ no vaccinated gets infected), and $0\leq\varepsilon \leq 1$.

Thus we get the following controlled SVIR model (see Fig. \ref{contr_svir_grafo}):
\begin{figure}[t]
\begin{center}
	\includegraphics[width=9cm, height=5cm]{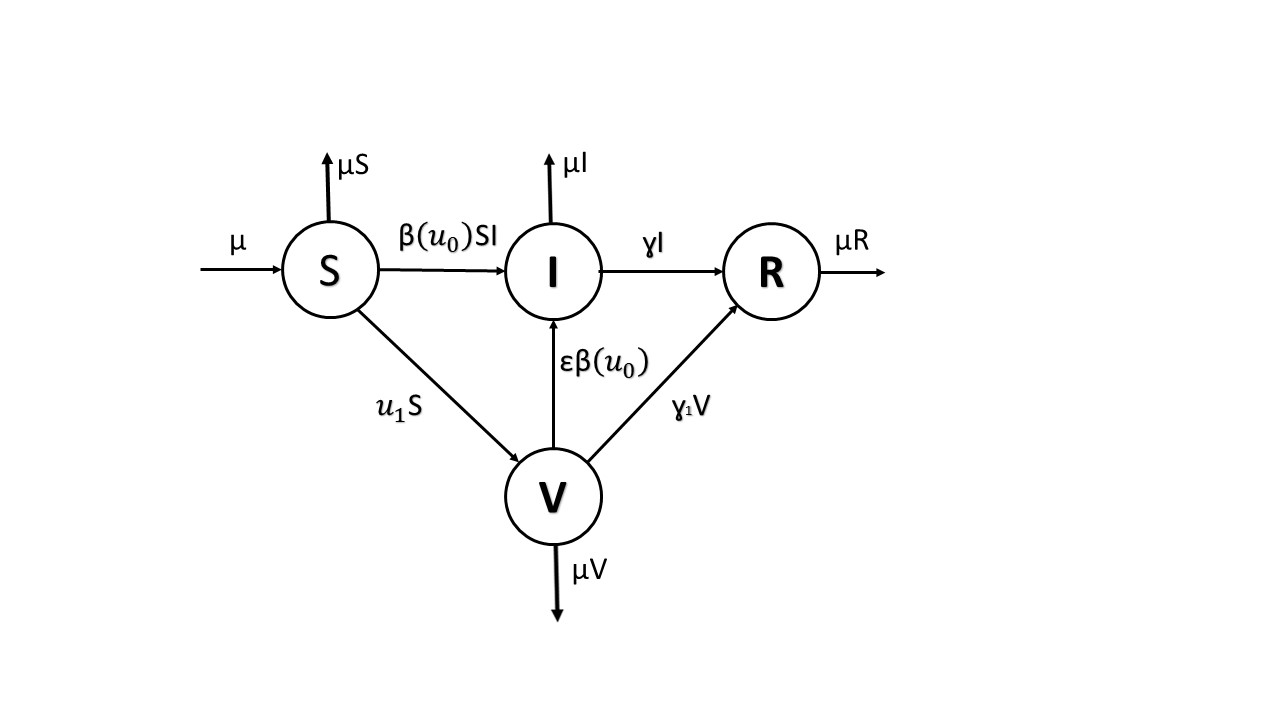}
\end{center}
\caption{The controlled SVIR model graph.}
\label{contr_svir_grafo}
\end{figure}

\begin{equation}
\label{sistema}
\left\{
\begin{array}{ll}
\displaystyle{\frac{dS}{dt}(t)} =-\beta (u_0(t))S(t)I(t) -u_1(t)S(t)+ \mu - \mu S(t) &S(0)=S_0
\\\\
\displaystyle{\frac{dV}{dt}}(t) =u_1(t)S(t)-\varepsilon \beta (u_0(t))V(t)I(t)-\gamma _1V(t)- \mu V(t) &V(0)=V_0
\\\\
\displaystyle{\frac{dI}{dt}}(t) =\beta (u_0(t))(t)S(t)I(t) +\varepsilon \beta (u_0(t))V(t)I(t)-\gamma I(t) - \mu I(t) &I(0)=I_0
\\\\
\displaystyle{\frac{dR}{dt}}(t)=\gamma _1V(t)+\gamma I(t)- \mu R(t) &R(0)=R_0
\end{array}
\right.
\end{equation}

where we assume that the initial data   $S_0, V_0, I_0, R_0 \in \mathbb{R}^+$, and $S_0+V_0+I_0+R_0 = 1$. The above assumptions are stated since the model (\ref{sistema}) represents human populations, and it can be shown that the solutions of the system are non--negative given non--negative initial values, see \cite{XTI}. As before, we immediately have from (\ref{sistema}) that $\frac {dN}{dt}(t) = 0$: hence $N(t)=N_0 \equiv 1$, for all $t\geq 0$.

Since the recovered people $R$ do not appear in the first three equations of (\ref{sistema}) or in the costs of the disease, we consider the following system
\begin{equation}
\label{sistemar}
\left\{
\begin{array}{ll}
\displaystyle{\frac{dS}{dt}}(t) =-\beta (u_0(t))(t)S(t)I(t) -u_1(t)S(t)+ \mu - \mu S(t) &S(0)=S_0
\\\\
\displaystyle{\frac{dV}{dt}}(t) =u_1(t)S(t)-\varepsilon \beta (u_0(t))V(t)I(t)-\gamma _1V(t)- \mu V(t) &V(0)=V_0
\\\\
\displaystyle{\frac{dI}{dt}}(t) =\beta (u_0(t))(t)S(t)I(t) +\varepsilon \beta (u_0(t))V(t)I(t)-\gamma I(t) - \mu I(t) &I(0)=I_0.
\end{array}
\right.
\end{equation}

To specify the control problem, we must define a functional to quantify the cost of spreading the disease. Specifically, we consider the costs due to infection and those attributable to vaccination, the latter assumed to be proportional to $u^2_1S$, the rate of individuals moving from $S$ to $V$ or, equivalently, the number of new vaccinated people per unit time. We categorize these expenses as comprising hospitalization costs for both inpatients, whether or not they require Intensive Care Unit (ICU) services, and logistical expenditures associated with the vaccination program, such as setting up and operating a vaccination hub, along with its medical staff, among other components. We finally include the cost of social restrictions in our framework, which we assume to be a function $c$ of the control $u_0$, such that $c$ is strictly increasing and convex in $u_0$, and that $c(0)=0$. This implies that, without control, the total costs of the disease spread are due to the infection and the vaccination. By assuming an additive structure for the cost functional, we separate the costs solely due to the disease from those due to the ``restrictions'' imposed on society. Parameters $c_1, c_2 \in \mathbb{R}^+$ represent the infection cost and the vaccination cost respectively.

Hence the objective function is given by $J:\cal {U}\rightarrow \mathbb{R}$ such that
\begin{equation}
\label{funzionale}
J(u_{0},u_{1}) =\int_{0}^{T} [c(u_0(t))+c_1I(t)+c_2u^2_1(t)S(t)]dt
\end{equation}

The aim is to find the best strategy $u^*\in U$ and the related state variables $S^*, V^*$ and $R^*$ which minimizes (\ref{funzionale}) 

$$
\min_{u\in \cal {U}} J(u)\qquad \mbox{subject to (\ref{sistemar})}.
$$
To prove the existence of such a strategy $u^*$ we refer to \cite {RT23} \cite{FR}, \cite {LW} and  \cite{KPS}.

\begin{theorem}	
Let $\beta (\cdot)$ be a linear decreasing function and  let $c(\cdot)$ be a strictly convex, twice continuous differentiable function, such that $c'>0$ and $c(0)=0$. 
	
Then an optimal solution $ u^*$  for problem (\ref{sistemar})--(\ref{funzionale}) exists, i.e. there exists an optimal control $u^*\in \cal {U}$ such that $J(u^*)=\min J(u)$.
\end{theorem} 

\proof 
First of all, notice that the right hand side functions of system (\ref{sistemar})  are Lipschitz continuous with respect to the state variables, hence Picard--Lindelof Theorem ensures that there exist solutions to (\ref{sistemar}). By definition, the set  $[0,\overline{u_0}]\times [0,\overline{u_1}]$ is compact and convex and system
(\ref{sistemar}) is linear in the control variable $u$, then the result follows applying Theorem 4.1 and Corollary 4.1 pp. 68--69 in \cite{FR}.
\qed

\begin{remark}
By choosing a continuous, convex function $C(u,I,S)$ on $[0,\overline{u_0}]\times [0,\overline{u_1}]$ in (\ref{funzionale}), we can obtain a similar result as in Corollary 4.1 pp. 68--69 of \cite{FR}. We use $C(u,I,S)= c(u_0(t))+c_1I(t)+c_2u^2_1(t)S(t)$ to separate the costs of social restrictions, infection and vaccination. This also allows us to solve the problem numerically.
\end{remark}

We use the control theory in \cite{FR} or \cite{LW} to solve the optimal control problem. We define the Hamiltonian function $H$ and the co-state variables $\lambda _1(t)$, $\lambda _2(t)$ and $\lambda _3(t)$. We drop the time dependence of the state variables $S, V, I$, the control variables $u_0, u_1$ and the co-state variables, unless stated otherwise. The Hamiltonian function of (\ref{sistemar})--(\ref{funzionale}) is given by

\begin{equation}
\label{hamiltoniana}
\begin{array}{l}
H(t,S,V,I,u_{0},u_{1},\lambda _1,\lambda _2,\lambda _3) 
=c(u_0(t))+c_1I+c_2u^2_1S+
\lambda _1[-\beta (u_0)SI -u_1S+ \mu - \mu S ]
\\\\
+\lambda _2[u_1S-\varepsilon \beta (u_0)VI-\gamma _1V- \mu V]
+\lambda _3[\beta (u_0)SI +\varepsilon \beta (u_0)VI-\gamma I - \mu I]
\end{array}
\end{equation}

\begin{theorem} \label{optimal_system}
Let $(S^*,V^* , I^*, u^*)$ be an optimal solution  for problem (\ref{sistemar})--(\ref{funzionale}), then there exist adjoint functions $\lambda _1,\lambda _2$  and $\lambda _3$ satisfying the following system of differential equations
\begin{equation} \label{costato}
\left\{
\begin{array}{l}
\lambda _1' =[\beta (u_0^*)I^* +u_1^*+ \mu]\lambda _1 - 
u_1^*\lambda _2-\beta (u_0^*)I^*\lambda _3-c_2{u^2_1}^*
\\\\
\lambda _2' =[\varepsilon \beta (u_0^*)I^* +\gamma _1+ \mu]\lambda _2 
-\varepsilon \beta (u_0^*)I^*\lambda _3
\\\\
\lambda _3' =\beta (u_0^*)S \lambda _1+\varepsilon \beta (u_0^*)V^*\lambda _2-
u_1^*\lambda _2-[\beta (u_0^*)S^*+\varepsilon \beta (u_0^*)V^*-\gamma -\mu]\lambda _3-c_1
\end{array}
\right.
\end{equation}
with the tranversality conditions on the co--states $\lambda _1$, $\lambda _1$ and $\lambda _3$ given by:
 $$
 \lambda _1(T)=0,\qquad \lambda _2(T)=0,\qquad  \lambda _3(T)=0.
$$
 The optimal restriction policy $u^*$ is such that
\begin{equation} \label{u*}
\displaystyle{ u^*(t)\in \mbox{argmin}_{u\in [0,1]\times [0,1]} H(t,S^*,V^*,I^*,u,\lambda _1,\lambda _2,\lambda _3) }.
\end{equation}
\end{theorem}

\proof
Let  $(S^*,V^* , I^*, u^*)$ be an optimal solution for problem (\ref{sistemar})--(\ref{funzionale}).  By Pontryagin's Maximum Principle the costate  variables $\lambda _1$, $\lambda _2$ and $\lambda _3$   satisfy system (\ref{costato}) whose equations are obtained evaluating the partial derivatives of the Hamiltonian function $H$ in (\ref{hamiltoniana}), with respect to the state variables $S,V,I$  

Suppose $(S^*,V^* , I^*, u^*)$ is an optimal solution of (\ref{sistemar})--(\ref{funzionale}). By Pontryagin's Maximum Principle, the co-state variables $\lambda _1$, $\lambda _2$ and $\lambda _3$ satisfy (\ref{costato}), which is derived from the partial derivatives of H in (\ref{hamiltoniana}) with respect to $S,V,I$.

\begin{equation}
\left\{
\begin{array}{l}

\lambda _1' =-\displaystyle{\frac{\partial H}{\partial S}}
\\\\
\lambda _2' =-\displaystyle{\frac{\partial H}{\partial V}}
\\\\
\lambda _3' =-\displaystyle{\frac{\partial H}{\partial I}}
\end{array}
\right.
\end{equation}
with the transversality conditions 
$\lambda _1(T)=\lambda _2(T)=\lambda _3(T)=0.$

The Hamiltonian function $H$, defined in (\ref{hamiltoniana}) is strictly convex with respect to the control variable $u$, hence the existence of a unique minimum follows, see \cite{WMV}, therefore

The strictly convexity, with respect to the control variable $u$, of the Hamiltonian function $H$, defined in (\ref{hamiltoniana}), ensures the existence of a unique minimum, see \cite{WMV}, hence
	
$$
u^*(t)\in \mbox{argmin}_{u\in [0,1]\times [0,1]} H(t,S^*,V^*,I^*,u,\lambda _1,\lambda _2,\lambda _3).
$$
\qed

\subsection{The solution of a class of optimal control problems}

We can make the result more specific by choosing a particular form for the transmission rate $\beta(u_0)$ and the cost function $c(u_0)$. For the former, we use a simple linear model:

\begin{equation} \label{lin_contr}
\beta(u_0) = \beta_0 (1- u_0), \ \ \ 0 \leq u_0 \leq 1.
\end{equation}

We assume that $\beta_0>0$ is the specific rate of how the disease spreads. In this case, we model the scenario when the maximum control (i.e. $u_0 \equiv 1$) completely stops the disease from diffusing.

For the social cost function $c(u_0)$, we choose the following specification for our practical application:
\begin{enumerate}
	\item $c_{quad}(u_0) = b {u_0}^2$, $ b >0$;
	\item $c_{exp}(u_0) = e^{k u_0}-1$, $k >0$;
\end{enumerate}

We prove the following result that gives a full description of the optimal controls.

\begin{proposition}
Let $\beta (u_0) = \beta_0 (1- u_0)$ and $c_{quad}(u_0)=b {u_0}^2$. Then the optimal control strategy $u_{quad}^*=(u_0^*,u_1^*)$ for  problem (\ref{sistemar})--(\ref{funzionale}) is given by
\begin{equation}	\label{u0*qdr}
\displaystyle{u_0^*(t)=min\left \{max \left [0,\frac{\beta_0I^*(t)[S^*(t)(\lambda _3(t)-\lambda _1(t))+ \varepsilon V^*(t)(\lambda _3(t)-\lambda _2(t))]}{2 b}\right ], \overline{u_0}\right \}}
\end{equation}
and
\begin{equation} \label{u1*qdr}
u_1^*(t)=min\left \{max\left [0,\frac{\lambda _1(t)-\lambda _2(t)}{2c_2}\right ], \overline{u_1}\right \}
\end{equation}
	
\end{proposition}

\proof
In this case the Hamiltonian function $H$ is defined as
\begin{equation} \label{hamiltoniana2}
\begin{array}{l}
H(t,S,V,I,u,\lambda _1,\lambda _2,\lambda _3) =bu_0^2+c_1I+c_2u_1^2 S+
\lambda _1[-\beta_0(1- u_0) S I -u_1S+ \mu - \mu S ]+
\\\\
+\lambda _2[u_1S-\varepsilon \beta_0 (1- u_0) VI-\gamma _1V- \mu V]
+\lambda _3[\beta_0 (1- u_0) SI +\varepsilon \beta_0 (1- u_0) VI-\gamma I - \mu I]
\end{array}
\end{equation}
then, imposing first order conditions to minimize the Hamiltonian $H$ at $S^*,I^*,V^*$
\begin{equation}
\begin{array}{l}
\displaystyle{\frac{\partial H}{\partial u_0}}=2bu_0 +I^*[\beta_0 S^*(\lambda _1-\lambda _3)+
\varepsilon \beta_0 V^*(\lambda _2-\lambda _3)]=0
\\\\
\displaystyle{\frac{\partial H}{\partial u_1}}=2c_2u_1S-\lambda_1S+\lambda_2S=0
\end{array}
\end{equation}
we derive the optimal  restriction policy $u_{quad}^*=(u_0^*,u_1^*)$.\qed

\vspace{0.5cm}

\begin{corollary}
	Let $\beta (u_0) = \beta_0 (1- u_0)$ and $c_{quad}(u_0)=b {u_0}^2$.
	
	 Then the optimal control strategy $u_{quad}^*=(u_0^*,u_1^*)$ for  problem (\ref{sistemar})--(\ref{funzionale}) satisfies
	$$
	\lim\limits_{b\to +\infty}u_0^*(t)=0 \qquad \mbox{and} \qquad \lim\limits_{b\to 0}u_0^*(t)=\overline{u_0} 
	$$
	$$
	\lim\limits_{c_2\to +\infty}u_1^*(t)=0 \qquad \mbox{and} \qquad \lim\limits_{c_2\to 0}u_1^*(t)=\overline{u_1} 
	$$

\end{corollary}

\proof From the previous Proposition, the dynamical system has bounded solutions for every $u \in \mathcal{U}$. By taking the limit in (\ref{u0*qdr}) and (\ref{u1*qdr}), we get the result.

\vspace{0.5cm}

We can find the solution for the exponential case in a similar way as we did for the quadratic case above.
\begin{proposition}
	\label{propexp}
	Let $\beta (u_0) = \beta_0 (1- u_0)$ and  $c_{exp}(u_0)=e^{k u_0}-1$.  
	
	If $\lambda _3(t)>\max \{\lambda _1(t),\; \lambda _2(t)\}$, then the optimal control strategy $u_{exp}^*(t)=(u_0^*,u_1^*)$ for  problem (\ref{sistemar})--(\ref{funzionale}) is given by

		\begin{equation}
		\label{u0*exp}
		u_{0}^*(t)=
		min\left \{max \left [0,\frac{1}{k}\ln \frac{\beta_0 I^*(t)K(t)}{k}\right ], \overline{u_0}\right \}
		\end{equation}
		and
		\begin{equation}
		\label{u1*exp}
	u_1^*(t)=min\left \{max\left [0,\frac{\lambda _1(t)-\lambda _2(t)}{2c_2}\right ], \overline{u_1}\right \}
		\end{equation}
	
	where K(t) is defined as
	$\displaystyle{
		K(t)=S^*(t)(\lambda _3(t)-\lambda _1(t))+
		\varepsilon V^*(t)(\lambda _3(t)-\lambda _2(t))}
	$.
\end{proposition}

\vspace{0.5cm}

\begin{remark} \label{costate_rem}
According to the Pontryagin's Maximum Principle, the optimal control is a function of three Lagrange multipliers: $\lambda _1$, $\lambda _2$ and $\lambda _3$. These multipliers correspond to the marginal cost of susceptible population, the marginal cost of vaccinated population, and the marginal cost of infected population, respectively. Furthermore, the difference between $\lambda _3$ and $\lambda _1$ can be interpreted as the marginal cost of having an additional susceptible person infected, while the difference between $\lambda _3$ and $\lambda _2$ can be interpreted as the marginal cost of having an additional vaccinated person infected. Moreover,  $\lambda _1-\lambda _2$ represents the marginal cost of having an extra susceptible vaccinated person, thus a person who even if vaccinated is still possible to be infected.
\end{remark}

\newpage

\section{A numerical study}

In this section, we present the results of a simulation study aimed at investigating the impact of an optimal control strategy on disease dynamics and the associated costs. The objective of this study is to assess the effectiveness of various control measures in mitigating the spread of the disease and reducing the economic burden. We highlight that this section serves the purpose of presenting the outcomes of the proposed model without engaging in any empirical investigation. In particular, no attempt has been made to fit the SVIR model to actual data. It is important to note that while an empirical analysis of this nature would be significant, it exceeds the scope of the current paper and is therefore left for future research.
%
To analyze the optimal control strategies in alternative settings, we rely on the application of the Forward-Backward Sweep (FBS) algorithm. Such an algorithm serves as an effective approach for solving optimal control problems by iteratively refining the control function based on the state and costate variables, see \cite{LW}. It consists of two main steps. Firstly, the forward state equations (\ref{sistemar}) are solved by utilizing an ordinary differential equation (ODE) solver. This step calculates the state variables based on the current control function $u_n(\cdot)$. Next, the costate equations (\ref{costato}) are solved backward in time using the same ODE solver. These equations represent the adjoint variables that provide information about the sensitivity of the cost function with respect to the state variables (see Remark \ref{costate_rem}). Based on the optimality conditions, the control function is updated using the calculated state and costate variables. This process generates a new approximation of the state, costate, and control $u_{n+1}(\cdot)$. These steps are repeated iteratively until a convergence criterion is satisfied, i.e. when the algorithm has reached an acceptable approximation of the optimal control function. See McAsey et al. (2012) \cite{MMH} for a deep analysis of the convergence properties of the FBS method.

\medskip

In the initial stage, we set up the algorithm by establishing the temporal discretization and determining the termination criterion. To be specific, we selected a fixed number of time points $N$, uniformly distributed over the time interval $[0, T]$. The termination criterion was defined based on the non-decreasing behavior of the cost functional (\ref{funzionale}). Furthermore, we incorporated a technique of weighted averaging to update the solution iteratively. This involved combining the new solution $u_{new}(\cdot)$ and the previous one $u_n(\cdot)$ in such a way $u_{n+1}(\cdot)=u_{new}(\cdot) (1-c^n)+u_n(\cdot) c^n$. In order to fine tune the algorithm hyperparameters, we ran a set of preliminary experiments on our baseline model (see Table \ref{baseline_param} and the description below) by changing the starting solution (no controls or full controls) and the smoothing parameter $c \in \mathcal{C}$, where $\mathcal{C} =\{c_i : c_i = c_{i-1}+ \Delta c, c_0=0, i=1, \ldots, n+1\}$. In our investigation, we found that a weighting constant value $c=0.99$ provides the lowest minimum of the cost functional (\ref{funzionale}), and ensures the smoothing properties of the resulting solution. In contrast, no effect on the final optimal solution is observed concerning the starting point. The description of the implemented algorithm is given in the following:
\floatname{algorithm}{Forward-Backward Sweep Algorithm:}
\begin{algorithm}
\renewcommand{\thealgorithm}{}
\caption{Let $\underline{y}(t) = (S(t), I(t), R(t))$ and $\underline{u}(t) = (u_0(t), u_1(t))$. We denote here with $\underline{\dot{y}}(t) = F(\underline{y}(t), \underline{u}(t))$  and $\underline{\dot{\lambda}}(t) = G(\underline{\lambda}(t), \underline{y}(t), \underline{u}(t))$ the ODE systems (\ref{sistemar}) and (\ref{costato}), respectively. Furthermore, the optimality conditions are written as $\underline{u}^*(t) = H(\underline{y}^*(t),\underline{\lambda}^*(t))$ (see formulas (\ref{u0*qdr}), (\ref{u1*qdr}) and (\ref{u0*exp}), (\ref{u1*exp})).}\label{algo}
\begin{algorithmic}
\Require $\underline{u}^{(0)}(t)$ and $J_0 = J(\underline{u}^{(0)})$ \Comment{see formula (\ref{funzionale})}
\Repeat{$\ \ \ n\geq 0$}  
\State $\underline{\dot{y}}^{(n+1)}(t) = F(\underline{y}^{(n+1)}(t), \underline{u}^{(n)}(t))$, \ \ \ $\underline{y}^{(n+1)}(0) = y_0$
\State $\underline{\dot{\lambda}}^{(n+1)}(t) = G(\underline{\lambda}^{(n+1)}(t), \underline{y}^{(n+1)}(t), \underline{u}^{(n)}(t))$, \ \ \ $\underline{\lambda}^{(n+1)}(T) = 0$
\State $\underline{u}^{(new)}(t) = H(\underline{y}^{(n+1)}(t),\underline{\lambda}^{(n+1)}(t))$
\State $\underline{u}^{(n+1)}(t) = (1-c^n) \underline{u}^{(new)}(t) + c^n \underline{u}^{(n)}(t)$
\Until{$J_{n+1} > J_n$}
\end{algorithmic}
\end{algorithm}
\newpage

Now we present results for different instances of the controlled SVIR model, in order to show the main features of our modeling framework. The first model we consider, named baseline model, is an extension of the one considered in \cite{RT23}.

\subsection{Numerical results: the baseline model, $R_0^C < 1$.}
To begin, we establish a set of fundamental epidemiological parameters for the simulations. The following parameters define the baseline model. We take our time unit to be a day. In particular, we set $\beta_0 = 0.22$, $\gamma = 0.0795$ and $\gamma_1 = 0.0714$ (implying $1/\gamma \approx 12.6$ and $1 / \gamma_1 \approx 14$ days respectively). The value of $\varepsilon$ can be set by using an
estimate of the vaccine effectiveness, $VE$, which is defined as the percentage reduction in risk of disease among vaccinated persons relative to unvaccinated persons, implying $\varepsilon \equiv (1-VE)$. In our experiments, we used the value of $VE$ as estimated for the three available vaccines for Covid19 (see  \cite{A2022} (Table 3)), implying $\varepsilon = 0.078$. 

As introduced in Section \ref{sect_control}, the cost functional (\ref{funzionale}) is given by the sum of three terms, each related to a specific aspect of the problem: the cumulative ``social cost'' $J_{SC}(u) = \int_0^T c(u_0(t)) dt$, ``infection cost'' $J_{IC}(u) = \int_0^T c_1 I(t) dt$, and ``vaccination cost'' $J_{VC}(u) = \int_0^T c_2 u_1^2(t) S(t) dt$. We chose the corresponding weights by normalizing with respect to the infection cost in such a way $c_1=1$ and $c_2=0.02$\footnote{These values were obtained by using data from an empirical investigation on the costs of the Covid19 disease outlined in Marcellusi et al. \cite{MFS22}, as described in Ramponi and Tessitore \cite{RT23}.}.

In order to get comparable results among the cost models, we set the parameters $b=0.04$, and $k=0.03922$ respectively, so that, in the given time period, the value of the social cost with the maximum control $J_{SC}(\bar u_0)$ is about the same. 
In particular, as a result of a set of preliminary numerical experiments, as fully described in Ramponi and Tessitore (2023) \cite{RT23} in the case of social control only, we found that by varying the social cost function parameter it is possible to identify a regime in which the full-control strategy produces lower costs than the no-control strategy. As the parameter increases, the full-control strategy becomes increasingly costly, and at the same time, the optimal strategy "converges" to the no-control strategy. 

\begin{table}[h] 
\small
\centering
\begin{tabular}{cclccl}
\hline
\multicolumn{3}{@{}c@{}}{Epidemiological/dynamical parameters} & \multicolumn{3}{@{}c@{}}{Economic parameters} \\ \hline
Parameter & Value & Description & Parameter & Value & Description \\ \hline \hline
$\beta_0$ & 0.220 &  transmission rate & $b$ & $0.4$ & quad. cost function \\
$\gamma$     & 0.0795 & recovery rate from infected & $k$ & $0.0392$ & exp cost function \\
$\gamma_1$ & 0.0714 &  full immunization rate & $c_1$ & $1$ & infection cost\\
$\epsilon$ & 0.078 & vaccine ineffectiveness & $c_2$ & $0.02$ & vaccination cost\\
$\bar u_0$ & 1 & maximum social control & & & \\
$\bar u_1$ & 0.006 & maximum vaccination rate & & & \\ 
$\mu$ & 2.5e-05 & daily death rate & & & \\ \hline
\end{tabular}
\caption{Parameters for the baseline model.} 
\label{baseline_param}
\end{table}

\medskip

In the framework of our baseline model, we consider three benchmark scenarios to compare with the optimal control: in the first one the disease is not controlled either by social restrictions or through vaccination (no-control/no-vax). This is the case where $u_0(\cdot)=u_1(\cdot)\equiv 0$. In the second scenario we assume an uncontrolled SVIR model: $u_0(\cdot)\equiv 0$ and a vaccination campaign at the highest possible rate $\bar u_1$, $u_1(\cdot)\equiv \bar u_1$. The third scenario is built by assuming a full social control, i.e. $u_0(\cdot) \equiv 1$, and $u_1(\cdot)\equiv \bar u_1$, as before. In particular, we assume that the health-care system is able to vaccinate about 90\% of the population in one year, that is $u_1(\cdot) \equiv \bar{u}_1 = 0.006$ (this means that $e^{-\bar{u}_1} \approx 0.1$). Regarding the starting conditions of the dynamic model, it is assumed that the infectious disease has been spreading in a population for a given period before containment measures are implemented and a vaccine has become available. Hence we set $I_0=0.04$, $V_0=0$, $R_0=0.12$ and $S_0=1-I_0-V_0-R_0$. 

\begin{figure}[h]
\hspace{-1.5cm}
\includegraphics[width=16cm, height=8cm]{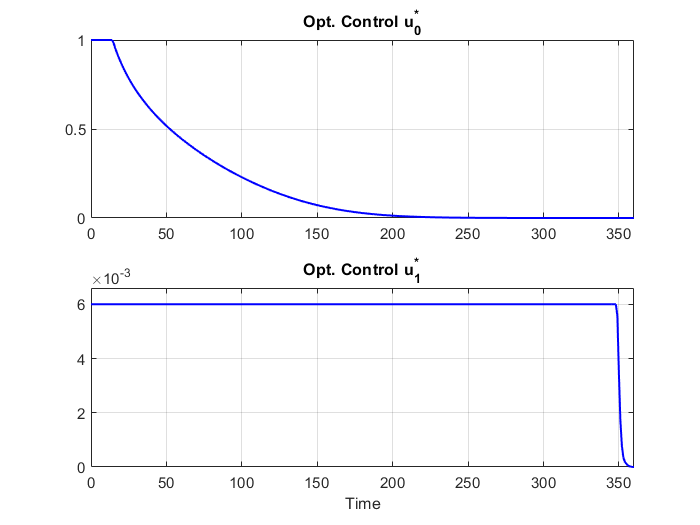}
\caption{Optimal controls $u_0^*, u_1^*$ for the baseline model, quadratic social cost function.}
\label{Fig_Optcontrols}
\end{figure}

\begin{figure}[h]
\hspace{-1.5cm}
\includegraphics[width=16cm, height=8cm]{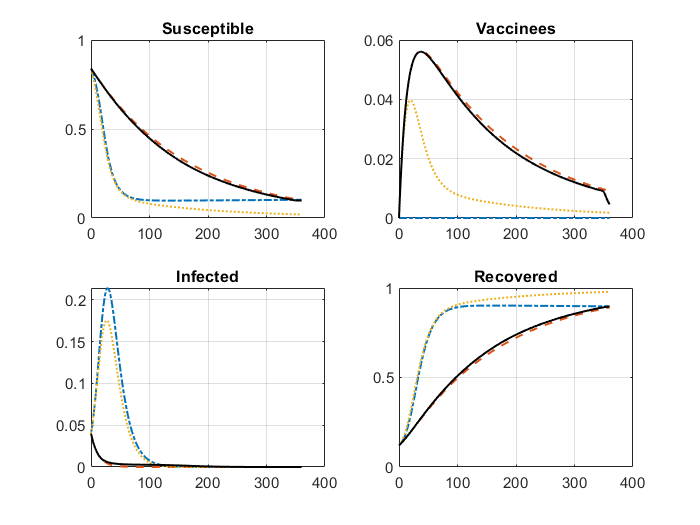}
\caption{Compartmental dynamics corresponding to the uncontrolled/no-vax system ($u_0 = u_1 \equiv 0$, dash-dotted blue line), uncontrolled SVIR ($u_0 \equiv 0, u_1 \equiv \bar u_1$, dotted yellow line), fully controlled system ($u_0 \equiv 1, u_1 \equiv \bar u_1$, dashed red line) and optimally controlled system ($u_0^*, u_1^*$, black line), quadratic social cost function.}
\label{Fig_Optcomparts}
\end{figure}

\begin{figure}[h]
\hspace{-1.5cm}
\includegraphics[width=16cm, height=8cm]{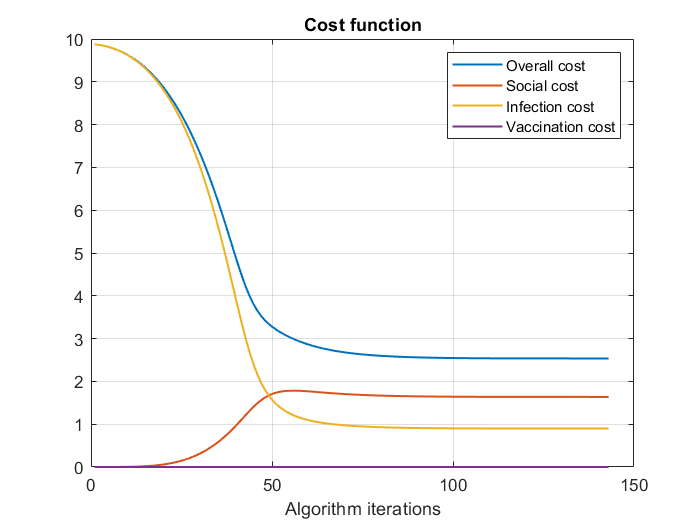}
\caption{Cost function values along the FBS algorithm iterations, quadratic social cost function.}
\label{Fig_FBS_conv}
\end{figure}

The results of these experiments are reported in Figures (\ref{Fig_Optcontrols}, \ref{Fig_Optcomparts}, \ref{Fig_FBS_conv}), 
for the quadratic social cost functions. As a matter of fact,  the qualitative behavior of the optimal solution for both the social cost functions is very similar. In particular, in both instances, the optimal controls suggest to maintain the maximal control of the transmission rate for few days and then to reduce it progressively; concurrently, the vaccination campaign should be implemented at the highest possible rate for almost all the considered period. In such a scenario, the the total cost is reduced considerably compared with the three benchmark strategies, as shown in Tables \ref{tabCosts_baseline_quad} and \ref{tabCosts_baseline_exp}, and the compartments dynamic is comparable with that of the benchmark with full control.  Specifically, the observed peak in the compartment of Infected individuals, under the absence of control measures, becomes entirely mitigated upon attaining, on the other hand, the same population percentage within the compartment of Recovered individuals, that is approximately $89.7 \%$, compared with the value $97.8 \%$ observed for the uncontrolled SVIR model. Moreover, the final percentage in the Susceptible compartment is around $10 \%$ in the case of the optimal strategy (as for the no-control/no-vax and full-control benchmarks), compared with $2 \%$ in the case of the uncontrolled SVIR model.

\begin{table}[h] 
\centering
\begin{tabular}{lc||ccc}
\hline
Contr. strategy & $J(u)$ &Social cost & Infection cost & Vaccination cost \\ \hline \hline
$u_0\equiv 0, u_1 \equiv 0$ &  $9.8731$  &    $0   (0\%)$ & $9.8731 (100\%)$  &  $0 (0\%)$ \\ 
$u_0\equiv 0, u_1 \equiv \bar u_1$ & $8.1198$ & $0 (0\%)$ & $8.1198 (99.99\%)$ & $0.00001 (0.0003\%)$ \\
$u_0\equiv 1, u_1 \equiv \bar u_1$ & $14.9033$    & $14.4000  (96.6226\%)$ & $0.5033  (3.3768\%)$ & $0.0001  (0.0006\%)$ \\
$u_0^*, u_1^*$ & $2.5362$   & $1.6369 (64.5393\%)$ &   $0.8993 (35.4573\%)$ &   $0.0001 (0.0034\%)$ \\
\hline
\end{tabular}
\caption{Total costs $J(u)$ of the strategies and the corresponding social, infection and vaccination cost for the quadratic social cost function.  In parentheses the percentage value w.r.t. the total cost.} 
\label{tabCosts_baseline_quad}
\end{table}

\begin{table}[h] 
\centering
\begin{tabular}{lc||ccc}
\hline
Contr. strategy & $J(u)$ &Social cost & Infection cost & Vaccination cost \\ \hline \hline
$u_0\equiv 0, u_1 \equiv 0$ & $9.8732$  &    $0   (0\%)$ & $9.8732 (100\%)$  &  $0 (0\%)$  \\ 
$u_0\equiv 0, u_1 \equiv \bar u_1$ & $8.1197$ & $0 (0\%)$ & $8.1197 (99.99\%)$ & $0.00001 (0.0003\%)$ \\
$u_0\equiv 1, u_1 \equiv \bar u_1$ & $14.9031$    & $14.3997 (96.6225\%)$ &   $0.5033 (3.3769\%)$ &  $0.0001 (0.0006\%)$ \\
$u_0^*, u_1^*$ & $3.5498$   & $2.6591 (74.9073\%)$ &   $0.8907 (25.0903\%)$ &   $0.0001 (0.0024\%)$ \\
\hline
\end{tabular}
\caption{Total costs $J(u)$ of the strategies and the corresponding social, infection and vaccination cost for the exponential social cost function.  In parentheses the percentage value w.r.t. the total cost.} 
\label{tabCosts_baseline_exp}
\end{table}

%
%

\paragraph{Sensitivity analysis.} The balance between the cost factors that determine the function $J$ affects the structure of the optimal solution. In particular, we analyzed the impact of the constraint $\bar u_1$ in relation to the cost of vaccination $c_2$. As expected when looking at the structure of the optimal control $u_1^*(\cdot)$, as the value of $c_2$ decreases, and $c_1$ and the parameter $b$ for the social cost function are fixed, the optimal control is "squeezed" toward the $\bar u_1$ constraint. 
\begin{figure}[h]
\hspace{-1.5cm}
\includegraphics[width=16cm, height=8cm]{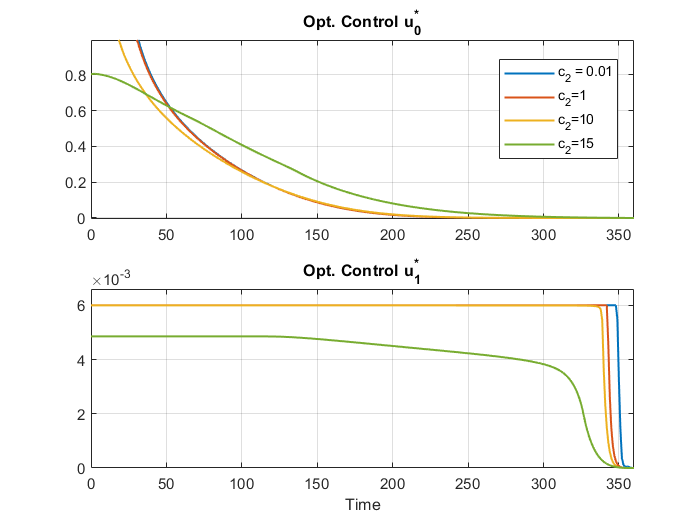}
\caption{Optimal strategies for different values of the vaccination cost $c_2$.}
\label{sens_c2}
\end{figure}

\subsection{Numerical results: the case $R_0^C >1$.}

\begin{figure}[h]
\hspace{-1.5cm}
\includegraphics[width=16cm, height=8cm]{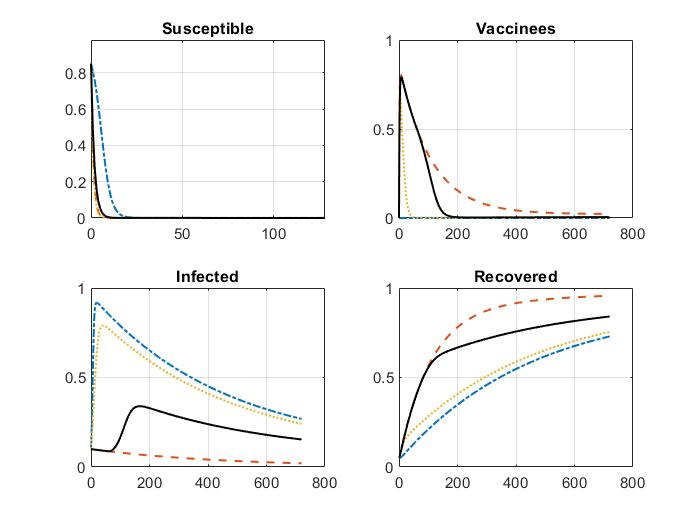}
\caption{Compartmental dynamics corresponding to the uncontrolled/no-vax system ($u_0 = u_1 \equiv 0$, dash-dotted blue line), uncontrolled SVIR ($u_0 \equiv 0, u_1 \equiv \bar u_1$, dotted yellow line), fully controlled system ($u_0 \equiv 1, u_1 \equiv \bar u_1$, dashed red line) and optimally controlled system ($u_0^*, u_1^*$, black line), endemic equilibrium model.}
\label{Fig_Optcomparts_endemic}
\end{figure}

In this section, we consider a set of epidemiological parameters implying a value $R_0^C >1$. As discussed in Section \ref{basic_svir}, such a condition implies the existence of an asymptotic equilibrium $(S^+,V^+,I^+,R^+)$ of the uncontrolled dynamical system for which $I^+ >0$. The parameters chosen in this scenario are reported in Table (\ref{param2}), and they imply a value $R_0^C =1.6261$.  We set as before $\bar u_0=1$, and an (unrealistically) high value for $\bar u_1=0.8$ so that the optimizer could determine an optimal vaccination strategy unaffected by this constraint. Time $T$ was set to $720$ days. 
\begin{table}[h] 
\small
\centering
\begin{tabular}{cclccl}
\hline
\multicolumn{3}{@{}c@{}}{Epidemiological/dynamical parameters} & \multicolumn{3}{@{}c@{}}{Economic parameters} \\ \hline
Parameter & Value & Description & Parameter & Value & Description \\ \hline \hline
$\beta_0$ & 0.4 &  transmission rate & $b$ & $0.12$ & quad. cost function \\
$\gamma$     & 0.002 & recovery rate from infected &  $c_1$ & $0.1$ & infection cost\\
$\gamma_1$ & 0.009 &  full immunization rate &  $c_2$ & $1$ & vaccination cost\\
$\epsilon$ & 0.4 & vaccine ineffectiveness &  & & \\
$\mu$ & 2.0e-04 & daily death rate & & & \\ \hline
\end{tabular}
\caption{Parameters for the endemic equilibrium model.} 
\label{param2}
\end{table}
Also, in this example, we can compare the behavior of the optimal system with benchmark models. In particular, we can observe how the susceptible compartment is quickly "emptied" due to a very high vaccination rate. On the other hand, the optimal strategy allows us to reduce the number of infected while increasing the Recovered considerably, see Figure \ref{Fig_Optcomparts_endemic}. In particular, the final values of Infected and Recovered compartments for the SVIR model is $24.14 \%$ and $75.47\%$, respectively, being instead $15.41\%$ and $84.06\%$ for the optimized model. At the same time, costs decrease significantly, as shown in Table \ref{tabCosts_endemic}. 
\begin{table}[h] 
\centering
\begin{tabular}{lc||ccc}
\hline
Contr. strategy & $J(u)$ &Social cost & Infection cost & Vaccination cost \\ \hline \hline
$u_0\equiv 0, u_1 \equiv 0$ &  $37.3778$  &    $0   (0\%)$ & $37.3778 (100\%)$  &  $0 (0\%)$ \\ 
$u_0\equiv 0, u_1 \equiv \bar u_1$ & $34.8572$ & $0 (0\%)$ & $33.3195 (95.5886\%)$ & $1.5377 (4.4114\%)$ \\
$u_0\equiv 1, u_1 \equiv \bar u_1$ & $91.6746$    & $86.4000  (94.2464\%)$ & $3.6129  (3.9411\%)$ & $1.6617  (1.8126\%)$ \\
$u_0^*, u_1^*$ & $26.7954$   & $9.8912 (36.9137\%)$ &   $15.7871 (58.9170\%)$ &   $1.1172 (4.1694\%)$ \\
\hline
\end{tabular}
\caption{Total costs $J(u)$ of the strategies and the corresponding social, infection and vaccination cost for the endemic equilibrium model.  In parentheses the percentage value w.r.t. the total cost.} 
\label{tabCosts_endemic}
\end{table}
\begin{figure}[h]
\hspace{-1.5cm}
\includegraphics[width=16cm, height=8cm]{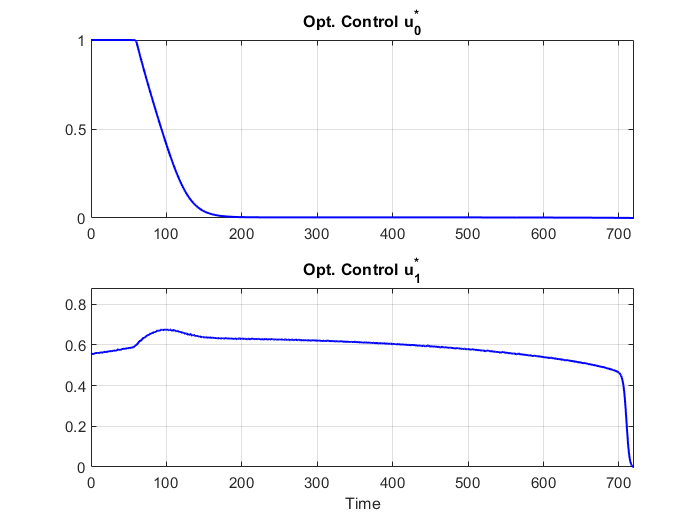}
\caption{Optimal controls $u_0^*, u_1^*$ for the endemic equilibrium model.}
\label{Fig_Optcontrols_endemic}
\end{figure}

The optimal strategy in this example (see Figure  \ref{Fig_Optcontrols_endemic}) sets the maximum social control in the initial period: this control gradually descends and simultaneously increases the optimal vaccination rate, which remains considerably high until it rapidly declines in the final period.

\section{Conclusions}

In this paper we considered a controlled SVIR compartmental model used to characterize the dynamics of infectious diseases in order to contribute to the development of a mathematical framework to analyze effective strategies for disease management to help policymakers and public health authorities in their decision-making processes. The model incorporates two distinct control mechanisms: i) social containment measures, encompassing a range of interventions such as lockdowns, curfews, school and university closures, and the cessation of commercial activities. These actions are effective in curbing disease spreading by lowering the transmission rate, but come at the expense of a social cost. ii) Vaccination campaign: The rate and efficiency of the vaccination campaign play a pivotal role in disease management. Vaccination not only reduces the spread of the disease but also incurs its own cost.

We explicitely considered the cost of the disease as the sum of three distinct terms: i) Social Cost, arising from the implementation of social containment measures, including economic and societal disruptions. ii) Infectious cost: Encompassing the toll exacted by the disease itself, in terms of morbidity, mortality, and healthcare burden. iii) Vaccination cost: Reflecting the expenses incurred in conducting the vaccination campaign. By using the Pontryagin Maximum Principle, 
under some condition on the functional form of these costs we were able to find the explicit structure of the optimal controls.

To assess the effectiveness of the proposed strategies, the optimally controlled system is simulated through the utilization of the Forward-Backward Sweep algorithm. This simulation approach facilitates a deeper understanding of the dynamics of infectious diseases and allows for the evaluation of the proposed control measures in practical scenarios.

In our numerical experiment we compared the system under the optimal strategy with three benchmark strategies: no controls /no vax, no social control /vax at maximum rate, and maximum social control / maximum vax rate. In our simulation, we observed that in the disease-free scenario ($R_C < 1$), the optimally controlled system significantly reduces the overall cost while maintaining the final values of each compartment almost on par with those of the fully controlled system. Conversely, in the endemic scenario ($R_C > 1$), the total cost is also reduced, but the final compartmental values fall between those achieved with the fully controlled system and those observed with other benchmark control strategies.

Finally, in future research, there is potential to explore various aspects of the SVIR compartmental model and its optimal control strategies. From an economic perspective, researchers may consider providing a more detailed description of the social cost function. This could involve investigating its relationship with financial and social indices, as well as refining the characterization of infectious costs, including distinctions between costs related to hospitalization with or without ICU care. From a dynamical point of view, future studies could expand the model by introducing additional compartments, such as "Quarantined" and "Dead." Incorporating these compartments would lead to significant changes in the dynamics of the compartmental ODE system, providing new insights into disease dynamics and control strategies. 
Last but not least, the calibration of this family of models to observed data, such as in Cerqueti et al. \cite{CRS23}, or Al\'os et al. \cite{AMS20} indeed represents a highly relevant and significant challenge.

\vspace{1cm}

\textbf{Author Contributions:} Both authors equally contributed  in every aspect of the \indent writing of this article. All authors have read and agreed to the published version of 
\indent the manuscript.



\textbf{Conflicts of Interests:} The authors declare no conflict of interest.

\end{document}